\documentclass[a4paper,11pt]{amsart}

\usepackage[utf8]{inputenc}
\usepackage{amsmath,amsfonts,amscd,amssymb,amsthm,epsfig,euscript}
\usepackage{mathrsfs}
\usepackage{amsxtra}
\usepackage{hyperref}
\hypersetup{colorlinks,citecolor=blue,urlcolor=blue,plainpages=false,hypertexnames=false}

\newtheorem{pro}{Proposition}

\theoremstyle{remark}

\numberwithin{equation}{section}
\usepackage[all]{xy}

\newcommand{\cal}{\mathcal}

\title{On $2k$-Hitchin's equations and Higgs bundles: a survey}

\author[Cardona]{S. A. H. Cardona}
\address{Conacyt Research Fellow--Instituto de Matem\'aticas, Universidad Nacional \indent Aut\'onoma de 
M\'exico, Leon 2 altos, Col. centro, 68000, Oaxaca, Mexico.}
\email{sholguin@im.unam.mx} 

\author[Garc\'ia-Compe\'an]{H. Garc\'ia-Compe\'an}
\address{Departamento de F\'isica, Centro de Investigaci\'on y de Estudios Avanzados \indent del IPN, P.O. Box. 14-740, C.P. 07000, Ciudad de M\'exico, Mexico.}
\email{compean@fis.cinvestav.mx} 

\author[Mart\'inez-Merino]{A. Mart\'inez-Merino}
\address{Conacyt Research Fellow--Facultad de Ciencias en F\'isica y Matem\'aticas, \indent Universidad Aut\'onoma de Chiapas, Carr. Emiliano Zapata Km. 08, Rancho \indent San Francisco, Col. Ter\'an, Tuxtla Guti\'errez, Mexico.}
\email{a.merino@fisica.ugto.mx; amartinezm@conacyt.mx}

\begin{document}

\maketitle 

\begin{abstract}
We study the $2k$-Hitchin equations introduced by Ward \cite{Ward 2} from the geometric viewpoint of Higgs bundles. After an introduction on Higgs bundles and $2k$-Hitchin's equations, we review some elementary facts on complex geometry and Yang-Mills theory. Then we study some properties of  holomorphic vector bundles and Higgs bundles and we review the Hermite-Yang-Mills equations together with two functionals related to such equations. Using some geometric tools we show that, as far as Higgs bundles is concern, $2k$-Hitchin's equations are reduced to a set of two equations. Finally, we introduce a functional closely related to $2k$-Hitchin's equations and we study some of its basic properties.\\


\end{abstract}


\section{Introduction}\label{Introduction}

This is a review article on Higgs bundles and a set of equations in mathematical physics, recently introduced by Ward \cite{Ward 2} and which are usually known as $2k$-Hitchin's equations. The purpose of the article is to begin a study of these equations using complex geometry. In order to do that, it is important to review some elementary facts on complex geometry and Yang-Mills theory; in particular, some properties of the Hodge operator, holomorphic vector bundles and Higgs bundles are crucial to develop a geometric approach of the $2k$-Hitchin equations. After studying these properties, we review the $2k$-Hitchin equations in the context of holomorphic vector bundles following the ideas of Kobayashi \cite{Kobayashi}. These can be seen as a set of four equations involving a metric in a holomorphic vector bundle $E$ and certain holomorphic form with values in the bundle of endomorphisms of $E$. From the above  geometric point of view, Higgs bundles seem to be a ``natural frame'' to study the $2k$-Hitchin equations. In fact, for these bundles the equations can be further reduced to a set of two equations defined for hermitian metrics of the bundle. Since Higgs bundles and $2k$-Hitchin's equations are of certain interest in mathematical-physics and complex geometry, the authors hope that the content of the present survey will be of some interest for physicists as well as mathematicians working in Yang-Mills theory and complex geometry.\\ 

The article is organized as follows. In Section \ref{Introduction} we give an historical introduction to Hitchin's equations, Higgs bundles and $2k$-Hitchin's equations. This section is not intended to be rigorous or exhaustive, however we would like to give a general overview on these topics. For more details on Higgs bundles and $2k$-Hitchin's equations see the articles of Simpson \cite{Simpson, Simpson 2} and Ward \cite{Ward 2}, respectively. A reader familiarized with supersymmetric Yang-Mills theories and string theory can find an introduction to Higgs bundles in \cite{Wijnholt}. In Section \ref{Hodge op.} we review some aspects of the Hodge $*$ operator in complex geometry. More precisely, we include some definitions and properties of this operator following the notation of Kobayashi \cite{Kobayashi}. In particular, we show that $*$ defines a natural hermitian inner product in the space of forms of a compact K\"ahler manifold $X$, this result appears in the form of Proposition \ref{Proposition A}.  The Hodge $*$ operator can be naturally extended to Yang-Mills theory, which is the main topic of Section \ref{CGeo-HYM}. In such a section we review some elementary definitions concerning holomorphic vector bundles, e.g., the notions of degree, hermitian metric and Chern connection of a holomorphic vector bundle $E\longrightarrow X$, where $X$ is a compact K\"ahler manifold. In particular, we define the hermitian adjoint of a form on $X$ with values in the bundle of endomorphisms of $E$. Using this we show that $*$ defines an hermitian inner product on the space of forms on $X$ with coefficients in the bundle of endomorphisms of $E$; this result appears as Proposition \ref{Proposition A(EndE)}. In fact, this hermitian inner product is a fundamental notion in any geometric approach to Yang-Mills theory, for details see the pioneering work of Atiyah \cite{Atiyah79}. \\

In Section \ref{Higgs-HYM section} we review the basic definitions on Higgs bundles. In particular, we introduce the Hitchin-Simpson connection and curvature and we define the Hermite-Yang-Mills equations for such bundles. We also introduce a couple of functionals closely related to such equations; namely the full Yang-Mills-Higgs functional and the Kobayashi functional. As a consequence of a result in \cite{Cardona 7} we know that there exists a non trivial expression relating these functionals, this result appears as  Proposition \ref{Prop. FullYM vs Koba}. In Section \ref{Higgs-Ward section} we review the $2k$-Hitchin equations of Ward \cite{Ward 2} from the point of view of a holomorphic vector bundle $E\longrightarrow X$, with $X$ a compact K\"ahler manifold and we show that such equations can be regarded as a set of four equations whose variables are pairs $(h,\Phi)$, with $h$ an hermitian metric in $E$ and $\Phi$ an holomorphic form of type $(1,0)$ on $X$ with values in the bundle of endomorphisms of $E$. At this point, we will see that two of these equations can be ``formally'' solved if we consider the Chern connection and a Higgs field on $E$. Therefore, as far as Higgs bundles is concern, the $2k$-Hitchin equations can be reduced to a set of only two equations defined in the space of hermitian metrics in the Higgs bundle. The remaining equations have some resemblance with Seiberg-Witten equations, and hence we propose a natural functional ${\cal H}(h)$ associated to the $2k$-Hitchin equations on a Higgs bundle; we call such a functional the non-abelian Seiberg-Witten functional. Finally, we show that a solution of $2k$-Hitchin's equations for Higgs bundles is necessarily a minimum of the non-abelian Seiberg-Witten functional, this result appears in the form of Proposition \ref{Prop. 2k-sol. minimum of H}. For the benefit of the reader, we include a Section \ref{Appendix} containing some remarks on Hitchin's equations, for more details see the articles of Hitchin \cite{Hitchin} and Ward \cite{Ward}. 


\subsection{Hitchin's equations}
Yang-Mills theory attracted the attention of geometers since late seventies and there exists a huge literature in mathematical-physics as well as in complex geometry on this topic. In particular, from a geometric point of view there are some pioneering works on Yang-Mills theory by Atiyah, Bott, Singer, Drinfeld, Manin, Donaldson and Hitchin \cite{Atiyah79, Atiyah-Bott, Atiyah-Drinfeld-Hitchin-Manin, Atiyah-Singer-Hitchin, Atiyah-Ward, Donaldson-2, Hitchin 0, Hitchin, Manin} among other authors. In fact, Hitchin's equations arise for the first time in \cite{Hitchin} as a dimensional reduction to ${\mathbb R}^{2}$ of the self-dual Yang-Mills equations (SDYM) on ${\mathbb R}^{4}$. Since such equations have a conformal invariance, these can be studied using a geometric setting. More precisely, Hitchin considers principal $G$-bundles $P\longrightarrow X$ with $X$ a compact Riemann surface and the equations are defined for a {\it connection form} $A$ on $P$ and a $(1,0)$-form $\Phi$ on $X$ with values in the adjoint bundle ${\rm ad}P$. Using the above notation the Hitchin's equations are usually written as:
\begin{equation}
F_{A} + [\Phi\,,\Phi^{*}] = 0\,, \quad\quad\quad d''_{A}\Phi = 0\,, \label{Hitchin eqs. 0}
\end{equation}  
where $F_{A}$ is the curvature of $A$, the bracket $[\cdot\,,\cdot]$ is the usual commutator of forms with values in ${\rm ad}P$ and $d''_{A}$ is the anti-holomorphic covariant derivative induced by $A$ (see Section \ref{Appendix} for more details). From a physical point of view $A$ is interpreted as a {\it gauge potential}, $\Phi$ is the {\it Higgs field} obtained from the dimensional reduction procedure and $\Phi^{*}$ represents its usual hermitian conjugate (matricial adjoint). In fact, due to the origin of $\Phi$ Hitchin call it a Higgs field. There exists literature related to Hitchin's equations. that have been written in the last two decades. In particular, in complex geometry we want to mention articles of Dunajski and Hoegner, Mosna and Jardim, and Wentworth \cite{Dunajski-Hoegner, Mosna-Jardim, Wentworth}, and in mathematical-physics there exist important articles written by Kapustin, Witten and Ward \cite{Witten-Kapustin, Ward, Ward 2,  Witten 2}. The Hitchin equations can be written in different ways, a fact which is in essence due to the several forms in which the gauge potential and the Higgs field can be considered. For instance, in mathematical-physics literature it is common to consider such equations in a real form instead of a complex one. In a more geometric approach it is frequently used a complex form of these equations. We include in Section \ref{Appendix} some different ways in which the Hitchin equations are presented.\\

As we already mentioned, Hitchin introduces the equations (\ref{Hitchin eqs. 0}) for principal $G$-bundles $P\longrightarrow X$ with $X$ a compact Riemann surface. Using some geometric tools, he studies these equations in the special cases when $G=SO(3)$ and $SU(2)$. In particular, he shows that for $G=SU(2)$ the existence of (non-singular) solutions depends on topological conditions of the Riemann surface as well as some algebraic conditions (Mumford stability) of a pair $(E,\Phi)$. Here the pair consists of a certain rank-2 holomorphic vector bundle $E$ related to $P$ and a section $\Phi$ of an associated bundle to $E$. For more geometric details the reader can be see the pioneering work \cite{Hitchin}. In general, the existence of solutions of Hitchin equations depends on topological properties of $G$ and $X$. For instance, for $X={\mathbb R}^{2}$ Mosna and Jardim \cite{Mosna-Jardim} found (non-singular) solutions when $G=SO(2,1)$ and more recently Ward \cite{Ward} studied (singular) solutions for $G=SU(2)$.

 
\subsection{The $2k$-Hitchin equations}
There exists a lot of interest on Hitchin's equations (\ref{Hitchin eqs. 0}) in mathematical physics. For instance, such equations play an important role in the celebrated works of Witten and Kapustin \cite{Witten-Kapustin, Witten 2} on the geometric Langlands programme. Additionally, there exists a generalization of Hitchin's equations introduced by Ward \cite{Ward 2}. To be precise, Ward introduces, for any $k$-complex manifold $X$ and Lie group $G$, the equations  
\begin{equation}
D_{A}\Phi = 0\,, \quad\quad F^{1,1}_{A} + \frac{1}{4}[\Phi\wedge\Phi^{*}] = 0\,, \quad\quad [\Phi\wedge\Phi] = 0\,, \quad\quad F^{2,0}_{A} = 0\,, \label{2k-Hitchin eqs.}
\end{equation}
where $A$ and $\Phi$ are the {\it gauge} and {\it Higgs fields}, respectively. Here the bracket $[\cdot]$ is a common notation in physics for the commutator of forms with values in the Lie algebra ${\mathfrak g}$ of $G$. In other words 
\begin{equation}
[\Phi\wedge\Phi^{*}] = [\Phi\,,\Phi^{*}]\,, \quad\quad\quad [\Phi\wedge\Phi] = 2\,\Phi\wedge\Phi\,. \nonumber
\end{equation}
The study of (\ref{2k-Hitchin eqs.}) from the geometric viewpoint of Higgs bundles is the main purpose of the present work. Ward called (\ref{2k-Hitchin eqs.}) the $2k$-Hitchin equations and considered these as a generalization of (\ref{Hitchin eqs. 0}) to higher dimensional complex manifolds. Now, from a geometric setting, the equations (\ref{2k-Hitchin eqs.}) can be defined for a holomorphic vector bundle $E$ over $X$, where  $D_{A}$ is considered as a {\it connection} on $E$ with connection form $A$ and $\Phi$ is a $(1,0)$-form on $X$ with coefficients in the bundle ${\rm End}E$. As we will see in Section \ref{Higgs-Ward section}, if the form $\Phi$ is holomorphic it can be interpreted as a {\it Higgs field} in the Higgs bundle setting. By imposing some algebraic conditions on the Higgs field $\Phi$, Ward \cite{Ward 2} was able to find some explicit (non-singular) solutions to (\ref{2k-Hitchin eqs.}) when $X={\mathbb C}^{2}$ and the gauge group $G$ involved in $A$ and $\Phi$ is $SU(2)$. Ward shows that for $k=2$, the equations (\ref{2k-Hitchin eqs.}) are related to another set of equations -- commonly known as the {\it non-abelian Seiberg-Witten equations} -- and for which Dunajski and Hoegner \cite{Dunajski-Hoegner} found solutions using a generalized t'Hooft ansatz.  As it was noticed by Ward, in the lowest dimensional case, i.e., for $k=1$, the equations (\ref{2k-Hitchin eqs.}) are reduced to (\ref{Hitchin eqs. 0}). In fact, in such a case $\Phi\wedge\Phi=0$ and $F^{2,0}_{A}=0$ are satisfied due to dimensional restrictions and also $F^{1,1}_{A}=F_{A}$ and $D_{A}\Phi=d''_{A}\Phi$. At this point (and after rescaling $\Phi$) we can write the remaining two equations in (\ref{2k-Hitchin eqs.}) as (\ref{Hitchin eqs. 0}).


\subsection{Higgs bundles}
In complex geometry there exists another generalization of (\ref{Hitchin eqs. 0}). In fact, following the ideas of Hitchin and some geometric tools in Yang-Mills theory, Simpson \cite{Simpson, Simpson 2} introduces the notion of a {\it Higgs bundle} as a generalization of the pairs defined first by Hitchin and closely related to the equations (\ref{Hitchin eqs. 0}). Roughly speaking, Simpson defines a Higgs bundle as a pair ${\mathfrak E}=(E,\Phi)$ consisting of a rank-$r$ holomorphic vector bundle $E$ over an $n$-dimensional K\"ahler manifold $X$ and a certain holomorphic $(1,0)$-form $\Phi$ on $X$ with coefficients in the bundle ${\rm End}E$, called the Higgs field of the Higgs bundle. We can consider hermitian metrics in Higgs bundles. In fact, an hermitian metric in a Higgs bundle ${\mathfrak E}$ is by definition an hermitian metric $h$ in the corresponding holomorphic bundle $E$, i.e., following Kobayashi \cite{Kobayashi} it is a $C^{\infty}$-field of hermitian metrics in the fibers of $E$. A Higgs bundle with a fixed hermitian metric is sometimes called an {\it hermitian Higgs bundle} \cite{Cardona 3, Cardona 7}. A standard result in complex geometry \cite{Kobayashi, Lubke, Siu} asserts that any hermitian metric $h$ in $E$ defines a Chern connection $D_{h}$, using it and the Higgs field, Simpson \cite{Simpson} defines another connection ${\cal D}_{h}$ on the Higgs bundle ${\mathfrak E}$. From a geometric viewpoint, it is natural to consider the metric $h$ as a ``variable'' and look for solutions to the equation:
\begin{equation}
\hat{\cal K}_{h} = c\,h\,, \label{HYM-eq. 0}
\end{equation} 
where $\hat{\cal K}_{h}$ is, strictly speaking, the hermitian form associated to the mean curvature ${\cal K}_{h}$ of ${\cal D}_{h}$ and $c$ is certain constant depending on invariants of the bundle $E$ (see Section \ref{Higgs-HYM section} for details). Since ${\cal K}_{h}$ is a section of ${\rm End}E$, the equation (\ref{HYM-eq. 0}) is indeed a set of equations commonly called the {\it Hermite-Yang-Mills} (HYM) or {\it Hermite-Einstein} (HE) equations. This terminology is introduced by Kobayashi \cite{Kobayashi} and Uhlenbeck and Yau \cite{Uhlenbeck-Yau} for a ``similar'' equation for holomorphic vector bundles.\footnote{If $E\longrightarrow X$ is a holomorphic vector bundle over a compact K\"ahler manifold, the HE or HYM equation is ${\hat K}_{h}=c\,h$, where ${\hat K}_{h}$ is the mean curvature of the Chern connection $D_{h}$ of $E$; such equation can be considered as a generalization of the K\"ahler-Einstein equation. In fact, if $E=TX$ and $h=g$ and since $X$ is K\"ahler, then ${\hat K}_{h}={\rm Ric}_{g}$ and the HE equation becomes ${\rm Ric}_{g}=c\,g$. The same terminology is used for the Higgs bundle case, even when such equations are not the same as the classical ones in complex geometry.} In the lowest dimensional case, the HYM equations are reduced to (\ref{Hitchin eqs. 0}), hence these are also a generalization of Hitchin's equations. Since the original articles of Hitchin and Simpson \cite{Hitchin, Simpson}, Higgs bundles have played an important role in geometry and there exists a very vast literature on this topic in complex geometry and mathematical physics. In particular, we would like to emphasize that several properties of Higgs bundles can be naturally extended to a certain kind of principal bundles, usually known as Higgs $G$-bundles \cite{Bruzzo-Granha 0}. In general, principal bundles play an important role in gauge theory and hence this extension is of interest in areas of high energy physics like string theory and Yang-Mills theory. Now, a priori, it could be natural to review (\ref{2k-Hitchin eqs.}) in a more geometric setting and a Higgs bundle seems to be an appropriate geometric ``frame'' to do it. In fact (see Section \ref{Higgs-Ward section} for details), if we consider $2k$-Hitchin's equations (\ref{2k-Hitchin eqs.}) on a Higgs bundle ${\mathfrak E}=(E,\Phi)$, two of them are automatically satisfied. In fact, $\Phi\wedge\Phi=0$ and $F^{2,0}=0$ hold just from the definition of Higgs bundle and taking the Chern connection as the connection on the corresponding holomorphic vector bundle $E$; in such a way $2k$-Hitchin's equations (\ref{2k-Hitchin eqs.}) are reduced to a set of two equations.\\

{\bf Acknowledgements:} S. A. H. C. was partially supported by the CONACyT grant 256126. Part of this article was done during several visits of the authors at Centro de Investigaci\'on y de Estudios Avanzados del IPN (CINVESTAV) in Mexico city and Instituto de Matem\'aticas de la Universidad Nacional Aut\'onoma de M\'exico in Oaxaca (IMUNAM). S. A. H. C., respectively H. G-C. and A. M-M., wants to thank CINVESTAV and IMUNAM for the hospitality. 


\section{The Hodge operator in complex geometry}\label{Hodge op.}

The main purpose of this section is to fix the notation that we will use in the rest of the article as well as to review some basic properties of the Hodge $*$ operator in complex geometry. The main content of this section is standard and can be found (though written in different notation) in classical texts of complex geometry \cite{Demailly, Griffiths-Harris, Kobayashi}.\\

Throughout this article $X$ will be a compact K\"ahler manifold of complex dimension $n$ and K\"ahler form $\omega$. The volume of $X$ can be written in terms of the K\"ahler form as:  
\begin{equation}
{\rm vol\,}X = \int_{X}\omega^{n}/n! \label{vol X}
\end{equation}
and we denote by $\Omega_{X}^{1,0}$ the {\it holomorphic cotangent bundle} to $X$ and by $\Omega_{X}^{0,1}$ its complex conjugate bundle. If $\{z^{\alpha}\}_{\alpha=1}^{n}$ is a local coordinate system of $X$, then $\{dz^{\alpha}\}_{\alpha=1}^{n}$ and $\{d{\bar z}^{\beta}\}_{\beta=1}^{n}$ are local frame fields for $\Omega_{X}^{1,0}$ and $\Omega_{X}^{0,1}$ and we can write 
\begin{equation}
 \omega = i\sum g_{\alpha\bar\beta}\,dz^{\alpha}\wedge d{\bar z}^{\beta}\,. \label{g and omega z,bar z}
\end{equation}
We would like to emphasize that we are using the standard index notation in complex geometry \cite{Kobayashi}. Therefore, the ``bar'' above the index $\beta$ in \eqref{g and omega z,bar z} is just a way to recall that it is associated to the complex conjugate variable ${\bar z}^{\beta}$ and thus $\bar\beta$ should be considered as the same index $\beta$. Certainly, this could be cumbersome in some physics literature, where a label above an index does represent a different index, e.g., in supersymmetry 
$\dot{\alpha}$ and $\alpha$ are considered as different indices. Since this article is a geometric approach to Higgs bundles and $2k$-Hitchin's equations, it is more appropriate to adopt the index notation in complex geometry. \\

We denote by $\Omega_{X}^{p,q}$, with $0\le p,q\le n$, the bundle over $X$ obtained by taking ($p$ and $q$ times) wedge products of $\Omega_{X}^{1,0}$ and $\Omega_{X}^{0,1}$. Following \cite{Kobayashi} we denote by $A^{p,q}_{X}$ the space of all $C^{\infty}$-sections of $\Omega^{p,q}_{X}$, i.e., the elements in $A^{p,q}_{X}$ are $C^{\infty}$-forms of type $(p,q)$ on $X$. Now, many computations in complex geometry can be simplified if we use a multi-index notation and if we consider unitary local frame fields instead of holomorphic ones. Let $\{\theta^{\alpha}\}_{\alpha=1}^{n}$ be a unitary local frame field of $\Omega^{1,0}_{X}$, i.e., a frame of forms of type $(1,0)$ such that $g_{\alpha{\bar\beta}}=\delta_{\alpha{\bar\beta}}$. In such a frame 
\begin{equation}
 \omega = i\sum\theta^{\alpha}\wedge{\bar\theta}^{\alpha}\,. \label{g and omega}
\end{equation}
Additionally, let $A=(\alpha_{1}, ... ,\alpha_{p})$ and $B=(\beta_{1}, ... ,\beta_{q})$ be multi-indices ordered in a strictly increasing way, i.e., with $\alpha_{1}<\cdots<\alpha_{p}$ and $\beta_{1}<\cdots<\beta_{q}$ and denote 
\begin{equation}
\theta^{A} = \theta^{\alpha_{1}}\wedge\cdots\wedge\theta^{\alpha_{p}}\,, \quad\quad\quad\quad {\bar\theta}^{B} = {\bar\theta}^{\beta_{1}}\wedge\cdots\wedge{\bar\theta}^{\beta_{q}}\,. \label{A and B}
\end{equation} 
If $A$ is a multi-index, we denote by $A'$ the complementary multi-index of it, hence $A'=(\alpha_{p+1}, ... ,\alpha_{n})$ wth $\alpha_{p+1}<\cdots<\alpha_{n}$. We denote by $\sigma^{AA'}$ the sign of the permutation $AA'$. Notice that from elementary permutation theory we get
\begin{equation}
\sigma^{A'A}=(-1)^{p(n-p)}\sigma^{AA'}\,. \label{sA'A property}
\end{equation}
Using the above multi-index notation any $\phi\in A^{p,q}_{X}$ can be written as:
\begin{equation}
\phi = \sum\phi_{A{\bar B}}\,\theta^{A}\wedge{\bar\theta}^{B}\,,    \label{phi-comp.}
\end{equation}
where each $\phi_{A{\bar B}}$ is a ${\mathbb C}$-valued smooth function. Notice that since the multi-indices are ordered we do not need to include a constant term in the preceding formula. In fact, if we consider a multi-index notation where the indices are in general not ordered, we have to include a factor of $\frac{1}{p!q!}$ in (\ref{phi-comp.}). This is the way in which Kobayashi writes the forms in \cite{Kobayashi}. Following again the standard index notation in complex geometry, we put a ``bar'' in the multi-index $B$ of the components of $\phi$ and thus, strictly speaking, ${\bar B}$ is the same multi-index $B$ and not a different one.\\

By applying (\ref{phi-comp.}) to a $\psi\in A^{p,q}_{X}$ we have  
\begin{eqnarray*}
{\bar\psi} & = & \overline{\sum\psi_{A{\bar B}}\,\theta^{A}\wedge{\bar\theta}^{B}}\\
               & = & \sum\overline{\psi_{A{\bar B}}}\,{\bar\theta}^{A}\wedge\theta^{B}\\
               & = & \sum(-1)^{pq}\overline{\psi_{A{\bar B}}}\,\theta^{B}\wedge{\bar\theta}^{A}\\
               & = & \sum{\bar\psi}_{B{\bar A}}\,\theta^{B}\wedge{\bar\theta}^{A},
\end{eqnarray*}
and hence ${\bar\psi}$ is an element in $A^{q,p}_{X}$ whose components are given by 
\begin{equation}
{\bar\psi}_{B{\bar A}}=(-1)^{pq}\overline{\psi_{A{\bar B}}}\,.  \label{cc comp. psi}
\end{equation} 

At this point we introduce the Hodge $*$ operator, it is by definition the $A^{0,0}_{X}$-linear mapping given by 
\begin{equation}
*: A^{p,q}_{X}\longrightarrow A^{n-q,n-p}_{X}\,, \quad\quad\quad\quad  *\,(\theta^{A}\wedge{\bar\theta}^{B}) = i^{n}\epsilon^{AB}\theta^{B'}\wedge{\bar\theta}^{A'}\,, \label{* op.}
\end{equation}
where 
\begin{equation}
\varepsilon^{AB} = (-1)^{np + n(n+1)/2}\sigma^{AA'}\sigma^{BB'} = \pm\,1\,. \nonumber
\end{equation}
As we will see in a moment, the exponent in the definition of $\varepsilon^{AB}$ is introduced in order to simplify some properties of the $*$ operator. Notice that from the definition of $\varepsilon^{AB}$ and the identity (\ref{sA'A property}) we get
\begin{equation}
\varepsilon^{BA} = (-1)^{n(p+q)}\varepsilon^{AB}\,, \quad\quad\quad \varepsilon^{AB}\varepsilon^{B'A'} = (-1)^{n+p+q}\,. \label{eAB property}
\end{equation}
Using (\ref{phi-comp.}) and (\ref{* op.}) we have
\begin{eqnarray*}
*^{2}\phi &=& *\sum\phi_{A{\bar B}}\,i^{n}\varepsilon^{AB}\theta^{B'}\wedge{\bar\theta}^{A'}\\
              &=& \sum\phi_{A{\bar B}}\,i^{n}\varepsilon^{AB}*(\theta^{B'}\wedge{\bar\theta}^{A'})\\
              &=& (-1)^{n}\sum\phi_{A{\bar B}}\,\varepsilon^{AB}\varepsilon^{B'A'}\theta^{A}\wedge{\bar\theta}^{B}\,,
\end{eqnarray*}
and hence the second identity of (\ref{eAB property}) implies that for every $\phi\in A^{p,q}_{X}$ we have 
\begin{equation}
*^{2}\phi = (-1)^{p+q}\phi\,. \label{prop. * op.}
\end{equation}
Now, following \cite{Kobayashi} we define the {\it complex conjugate} ${\bar *}$ of $*$ as the operator given by 
\begin{equation}
{\bar *}: A^{p,q}_{X}\longrightarrow A^{n-p,n-q}_{X}\,, \quad\quad\quad\quad {\bar *}\,\phi =  *\,{\bar\phi}\,.      \label{bar * op.}
\end{equation}
The first identity of (\ref{eAB property}) implies that $*\,{\bar\phi} = \overline{*\,\phi}$ for every $\phi\in A^{p,q}_{X}$, we leave the proof of this as an exercise to the reader. Also a direct computation shows that
\begin{equation}
{\bar *}\,(\theta^{A}\wedge{\bar\theta}^{B}) = (-1)^{pq}i^{n}\epsilon^{BA}\theta^{A'}\wedge{\bar\theta}^{B'}\,. \label{prop. bar * op.}
\end{equation}
Using (\ref{prop. bar * op.}), the first identity of (\ref{eAB property}) and also that
\begin{equation}
{\bar\theta}^{B}\wedge\theta^{A'}=(-1)^{q(n-p)}\theta^{A'}\wedge{\bar\theta}^{B} \nonumber
\end{equation} 
we get
\begin{eqnarray*}
(\theta^{A}\wedge{\bar\theta^{B}})\wedge{\bar *}\,(\theta^{A}\wedge{\bar\theta}^{B})&=& (-1)^{pq}i^{n}\varepsilon^{BA}\theta^{A}\wedge{\bar\theta}^{B}\wedge\theta^{A'}\wedge{\bar\theta}^{B'}\\
                                                                                                                                    &=& (-1)^{np}i^{n}\varepsilon^{AB}\theta^{A}\wedge\theta^{A'}\wedge{\bar\theta}^{B}\wedge{\bar\theta}^{B'}\\
                                                                                                                                    &=& i^{n}(-1)^{n(n+1)/2}\sigma^{AA'}\sigma^{BB'}\theta^{A}\wedge\theta^{A'}\wedge{\bar\theta}^{B}\wedge{\bar\theta}^{B'},
\end{eqnarray*}
where in the last equality we have used the definition of $\varepsilon^{AB}$. Now, from elementary permutation theory we know that 
\begin{equation}
\sigma^{AA'}\theta^{A}\wedge\theta^{A'} = \theta^{1}\wedge\cdots\wedge\theta^{n}\,, \quad\quad\quad \sigma^{BB'}{\bar\theta}^{B}\wedge{\bar\theta}^{B'} = {\bar\theta}^{1}\wedge\cdots\wedge{\bar\theta}^{n}\,, \nonumber
\end{equation}
and replacing these identities in the expression above (and using again an argument of permutation theory) we obtain
\begin{equation}
(\theta^{A}\wedge{\bar\theta}^{B})\wedge{\bar *}\,(\theta^{A}\wedge{\bar\theta}^{B}) = i^{n}\,\theta^{1}\wedge{\bar\theta}^{1}\wedge\cdots\wedge\theta^{n}\wedge{\bar\theta}^{n}= \frac{\omega^{n}}{n!}\,. \label{prop. bar * op. 2}
\end{equation}
The identity (\ref{prop. bar * op. 2}) is a key property of the operator ${\bar *}$ and, as we will see, it will be useful not only in complex geometry but also in Yang-Mills theory. Notice that from (\ref{prop. bar * op.})      
and (\ref{prop. bar * op. 2}) we have that the $2n$-form 
\begin{equation}
(\theta^{A}\wedge{\bar\theta}^{B})\wedge{\bar *}\,(\theta^{C}\wedge{\bar\theta}^{D}) \nonumber
\end{equation} 
is different from zero if and only if $A=C$ and $B=D$. Hence 
\begin{eqnarray*}
\phi\wedge {\bar *}\,\psi &=& \sum\phi_{A{\bar B}}\,\overline{\psi_{C{\bar D}}}\,(\theta^{A}\wedge{\bar\theta}^{B})\wedge{\bar *}\,(\theta^{C}\wedge{\bar\theta}^{D})\\
                                    &=& \sum\phi_{A{\bar B}}\,\overline{\psi_{A{\bar B}}}\,(\theta^{A}\wedge{\bar\theta}^{B})\wedge{\bar *}\,(\theta^{A}\wedge{\bar\theta}^{B})\\
                                    &=& \sum\phi_{A{\bar B}}\,\overline{\psi_{A{\bar B}}}\,\frac{\omega^{n}}{n!}\,. 
\end{eqnarray*}
Since $\omega^{n}/n!$ is real the above expression implies that 
\begin{equation}
\overline{\phi\wedge {\bar *}\,\psi} = \psi\wedge {\bar *}\,\phi\,, \label{h.c. prop.} 
\end{equation}
and we have a local hermitian inner product $(\cdot\,,\cdot)$ on $A^{p,q}_{X}$ given by 
\begin{equation}
(\phi\,,\psi) = \sum\phi_{A{\bar B}}\,\overline{\psi_{A{\bar B}}}\,. \label{Local formula Ap,q}
\end{equation}
It is possible to rewrite (\ref{Local formula Ap,q}) in another way. In fact, by defining
\begin{equation}
 {\bar\psi}^{\gamma_{1}...\gamma_{p}{\bar\delta}_{1}...{\bar\delta}_{q}} = g^{\beta_{1}{\bar\delta}_{1}}\cdots g^{\beta_{q}{\bar\delta}_{q}} g^{\gamma_{1}{\bar\alpha}_{1}}\cdots g^{\gamma_{p}{\bar\alpha}_{p}}{\bar\psi}_{\beta_{1}...\beta_{q}{\bar\alpha}_{1}...\bar\alpha_{p}}  \label{bar-psi-contrav. comp.}
\end{equation}
and since we are considering unitary frames, $g^{\alpha{\bar\beta}}=\delta^{\alpha{\bar\beta}}$ and (\ref{bar-psi-contrav. comp.}) becomes
\begin{equation}
{\bar\psi}^{A{\bar B}}={\bar\psi}_{B{\bar A}}\,. \label{bar-psi-contrav.} 
\end{equation}
Therefore, from (\ref{cc comp. psi}) and (\ref{bar-psi-contrav.}) we can rewrite (\ref{Local formula Ap,q}) as:
\begin{equation}
(\phi\,,\psi)  = (-1)^{pq}\sum\phi_{A{\bar B}}\,{\bar\psi}^{A{\bar B}}\,. \label{Inner pto. forms E. Conv.}
\end{equation}
The formula (\ref{Inner pto. forms E. Conv.}) is written in a form which is more familiar with physics literature. In particular, notice that we could use the Einstein convention and rewrite  such a formula without using the summation symbol.\\

From (\ref{h.c. prop.}) and the above analysis we obtain the following classical result in complex geometry \cite{Kobayashi}.
\begin{pro}\label{Proposition A}
If $X$ is a compact K\"ahler manifold with K\"ahler form $\omega$, then 
\begin{equation}
\langle\phi\,,\psi\rangle = \int_{X}\phi\wedge {\bar *}\,\psi = \int_{X}(\phi\,,\psi)\,\frac{\omega^{n}}{n!} \,, \quad\quad\quad \forall\phi,\psi\in A^{p,q}_{X}\,\label{Inner pto. forms}
\end{equation}
gives a global hermitian inner product on $A^{p,q}_{X}$, where $(\phi\,,\psi)$ is locally given by (\ref{Local formula Ap,q}) or equivalently by (\ref{Inner pto. forms E. Conv.}). 
\end{pro} 

Since (\ref{Local formula Ap,q}) and (\ref{Inner pto. forms}) define local and global hermitian inner products, we have local and global norms given by the standard formulas 
\begin{equation}
\lvert\phi\lvert^{2} = (\phi\,,\phi)\,, \quad\quad\quad\quad  \|\phi\|^{2} = \langle\phi\,,\phi\rangle\,. \label{pw and L2 norms}
\end{equation}
It is important to note that the hermitian inner products $(\cdot\,,\cdot)$ and $\langle\cdot\,,\cdot\rangle$ depend on the K\"ahler metric $g=(g_{\alpha\bar\beta})$ of $X$. In fact, such a dependence appears in an implicit way in (\ref{* op.}) since the unitary local frame field $\{\theta^{\alpha}\}_{\alpha=1}^{n}$ depends on $g$. Now, in complex geometry the hermitian inner product defined by (\ref{Inner pto. forms}) is usually extended in a ``diagonal way'' to the space of all forms: it is zero when we evaluate $\langle\cdot\,,\cdot\rangle$ in forms of different type.


\section{Complex geometry and Yang-Mills theory}\label{CGeo-HYM}

In this section we study aspects of the differential geometry of complex vector bundles. In particular, we review some basic properties of holomorphic vector bundles over compact K\"ahler manifolds and we extend the Hodge $*$ operator to the space of $(p,q)$-forms on a compact K\"ahler manifold with values in the endomorphism bundle of a holomorphic vector bundle. The main purpose of this section is to fix a geometric framework in order to develop the following sections. \\

Let $E\longrightarrow X$ be a rank-$r$ holomorphic vector bundle, i.e., it is a complex vector bundle of rank $r$ with $E$ and $X$ complex manifolds and whose trivializations are biholomorphisms. Let $c_{1}(E)$ be its first Chern class, then the {\it degree} of $E$ is given by
\begin{equation}
{\rm deg}E = \int_{X}c_{1}(E)\wedge\omega^{n-1}. \label{deg E}
\end{equation}
Let us consider the bundle 
\begin{equation}
\Omega^{p,q}_{X}({\rm End}E) = \Omega^{p,q}_{X}\otimes{\rm End}E\,. \nonumber
\end{equation} 
It is the complex bundle whose sections are forms of type $(p,q)$ on $X$ with coefficients in ${\rm End}E$. The space of all these sections is usually denoted by $A^{p,q}_{X}({\rm End}E)$. As we mentioned in Section \ref{Introduction}, from a standard result in complex geometry \cite{Kobayashi, Lubke, Siu} we know that given an hermitian metric $h$ in the bundle\footnote{An hermitian metric $h$ in a complex vector bundle $E$ is a $C^{\infty}$-field of hermitian inner products in the fibers $E_{p}$ of $E$, i.e., for all $p\in X$ and $v,w\in E_{p}$, we have that $h(v,w)$ is linear in the first slot, $h(v,w) = \overline{h(w,v)}$ and $h(v,v)>0$ for $v\neq 0$ and also $h(s,s')$ is a $C^{\infty}$-function for all $C^{\infty}$-sections $s,s'$ of $E$. See \cite{Demailly} or \cite{Kobayashi} for details.} $E$, there exists a unique connection $D_{h}$ usually called the {\it hermitian} or {\it Chern connection}. It is the unique connection on $E$ compatible with the holomorphic structure and the hermitian metric $h$, i.e., decomposing $D_{h}=D_{h}' + D_{h}''$, it is the unique connection satisfying 
\begin{equation}
D_{h}'' = d'' = \sum{\bar\theta}^{\beta}\partial_{\bar\beta} \nonumber
\end{equation}
with $\partial_{\bar\beta}=\frac{\partial}{\partial{\bar\theta}^{\beta}}$, and
\begin{equation}
dh(s,s')=h(D_{h}s,s') + h(s,D_{h}s')\,, \quad\quad \forall s,s'\in\Gamma(E)\,. \nonumber
\end{equation}
The curvature of the Chern connection, denoted here by $F_{h}$, is usually called the Chern curvature. Hence $F_{h}=D_{h}\wedge D_{h}$ and it is known in complex geometry \cite{Kobayashi} that $F_{h}\in A^{1,1}_{X}({\rm End}E)$ and that  the equivalence class $c_{1}(E)$ can be represented by the form $\frac{i}{2\pi}{\rm tr\,}F_{h}$ in $A^{1,1}_{X}$.\\ 

Given an hermitian metric $h$ on $E$, we can define the {\it hermitian conjugate} (with respect to $h$) of any $\Psi\in A^{p,q}_{X}({\rm End}E)$, as the element ${\bar\Psi}_{h}\in A^{q,p}_{X}({\rm End}E)$ satisfying the condition:
\begin{equation}
h({\bar\Psi}_{h}s,s') = h(s,\Psi s')\,, \quad\quad \forall s,s'\in\Gamma(E)\,. \label{adjoint Psi}
\end{equation}
The hermitian conjugate ${\bar\Psi}_{h}$ can be considered as the adjoint of $\Psi$ with respect to $h$, sometimes it is denoted by $\Psi_{h}^{*}$ or just $\Psi^{*}$ \cite{Bruzzo-Granha, Cardona 3}. A straightforward computation using (\ref{adjoint Psi}) and the hemiticity of $h$ shows that the hermitian conjugate of $\bar\Psi_{h}$ is again $\Psi$ (we leave this as an exercise to the reader). Now, as we have seen in Section \ref{Hodge op.}, in complex geometry we have the operator
\begin{equation}
*: A^{p,q}_{X} \longrightarrow A^{n - q,n - p}_{X} \nonumber
\end{equation}
given by (\ref{* op.}). This map can be extended in the obvious way to a {\it Hodge operator} in Yang-Mills theory \cite{Atiyah79, Cardona 7}, and therefore we have in a natural way a map
\begin{equation}
*: A^{p,q}_{X}(\textrm{End}E) \longrightarrow A^{n - q,n - p}_{X}(\textrm{End}E)\,. \label{* op. EndE}
\end{equation}
Using the hermitian conjugation we can define the {\it hermitian conjugate} $\overline{*}_{h}$ of the Hodge $*$ operator (\ref{* op. EndE}), as the mapping given by the composition of (\ref{adjoint Psi}) and (\ref{* op. EndE}), i.e., 
\begin{equation}
\overline{*}_{h}: A^{p,q}(\textrm{End}\, E) \longrightarrow A^{n - p,n - q}({\rm End}E)\,, \quad\quad \overline{*}_{h}\,\Psi = *\left({\bar\Psi}_{h}\right). \label{adjoint *}
\end{equation}
The operator $\overline{*}_{h}$ can be seen as the extension of the usual ${\bar *}$ in complex geometry.\\

In Yang-Mills theory if $\Phi\in A^{p,q}_{X}({\rm End}E)$ and $\Psi\in A^{s,u}_{X}({\rm End}E)$ we can define a commutator 
$[\cdot\,,\cdot]$ of these forms as the element in $A^{p+s,q+u}_{X}({\rm End}E)$ given by the formula:
\begin{equation}
[\Phi\,,\Psi] = \Phi\wedge\Psi - (-1)^{(p+q)(s+u)}\Psi\wedge\Phi\,, \label{commutator}  
\end{equation} 
where the wedge product $\wedge$ is defined in the obvious way.\\

Using the preceding conventions we can write all these formulas in components. In particular, if $\Phi\in A^{p,q}_{X}({\rm End}E)$ we have 
\begin{equation}
\Phi = \sum\Phi_{A{\bar B}}\,\theta^{A}\wedge{\bar\theta}^{B}, \label{Phi comp.} 
\end{equation}
where the $\Phi_{A{\bar B}}$ are (locally defined) endomorphisms of $E$. These are the components of $\Phi$ with respect to the frame $\{\theta^{\alpha}\}_{\alpha=1}^{n}$. If $\{e_{i}\}_{i=1}^{r}$ is a local frame field of $E$ and $\{e^{j}\}_{j=1}^{r}$ is the dual frame field. Then we can write (\ref{Phi comp.}) in a more explicit way as: 
\begin{equation}
\Phi = \sum\Phi^{i}_{A{\bar B}j}\,e_{i}\otimes e^{j}\,\theta^{A}\wedge{\bar\theta}^{B}, \label{Phi comp 2.} 
\end{equation}
where the $\Phi^{i}_{A{\bar B}j}$ are the components of the endomorphisms $\Phi_{A{\bar B}}$ with respect to the frame field. Notice that these components define a ${\mathbb C}$-valued $r\times r$-matrix associated to $\Phi_{A{\bar B}}$. \\

Now, if $\Psi\in A^{s,u}_{X}({\rm End}E)$ we get 
\begin{eqnarray}
\Phi\wedge\Psi &=& \sum\Phi_{A{\bar B}}\Psi_{C{\bar D}}\,\theta^{A}\wedge{\bar\theta}^{B}\wedge\theta^{C}\wedge{\bar\theta}^{D} \nonumber \\ 
                        &=& (-1)^{qs} \sum\Phi_{A{\bar B}}\Psi_{C{\bar D}}\,\theta^{A}\wedge\theta^{C}\wedge{\bar\theta}^{B}\wedge{\bar\theta}^{D},   \label{Local wedge}
\end{eqnarray}
where the factor $(-1)^{qs}$ in the last equality appears as a consequence of the commutation of the terms with multi-indices $C$ and $B$. Notice that in a multi-index notation 
$\theta^{A}\wedge\theta^{C}=\pm\,\theta^{G}$, with $G=(\gamma_{1},...,\gamma_{p+s})$ the ordered set obtained from the $(p+s)$-tuplet $AC$. Hence, it is possible to write the above formula using an (ordered) multi-index notation with extra $\pm 1$ signs. In practice, it is not necessary to write such a general formula. Using this local expression and (\ref{commutator}) we have
\begin{eqnarray}
[\Phi\,,\Psi] &=&  \sum \, [\Phi_{A{\bar B}},\Psi_{C{\bar D}}]\, \theta^{A}\wedge \theta^{\bar B}\wedge \theta^{C}\wedge \theta^{\bar D} \nonumber \\
                &=& (-1)^{qs} \sum\, [\Phi_{A{\bar B}},\Psi_{C{\bar D}}]\, \theta^{A}\wedge\theta^{C}\wedge{\bar\theta}^{B}\wedge{\bar\theta}^{D}.  \label{Local commutator}
\end{eqnarray} 

According to (\ref{adjoint Psi}) we have
\begin{equation}
{\bar\Psi}_{h} = \sum{\bar\Psi}_{hB{\bar A}}\,\theta^{B}\wedge{\bar\theta}^{A} = \sum{\bar\Psi}^{i}_{hB{\bar A}j}\,e_{i}\otimes e^{j}\,\theta^{B}\wedge{\bar\theta}^{A}. \label{Psi-bar-h comp.}
\end{equation}
We can compute the components ${\bar\Psi}_{hB{\bar A}}$ of ${\bar\Psi}_{h}$ in terms of the components $\Psi_{A{\bar B}}$ of $\Psi$. In fact, if $\{e_{i}\}_{i=1}^{r}$ is again a local frame field of $E$, a straightforward  computation shows that (\ref{adjoint Psi}) implies 
\begin{equation}
{\bar\Psi}^{m}_{hB{\bar A}j} = (-1)^{pq}h_{j{\bar l}}\,\overline{\Psi^{l}_{A{\bar B}k}}\,h^{m{\bar k}},\label{bar-Psi vs. Psi gen.}
\end{equation}
where $h_{i{\bar k}}=h(e_{i},e_{k})$ and the components $h^{j{\bar k}}$ are defined by 
\begin{equation}
\sum h_{i{\bar k}}h^{j{\bar k}} = \delta_{i}^{\,j}\,.
\end{equation} 
In other words, as a matrix $(h^{j{\bar k}})$ is the inverse transpose of $(h_{i{\bar k}})$. In particular if $\{e_{i}\}_{i=1}^{r}$ is unitary, i.e., $h_{i{\bar k}}=\delta_{i{\bar k}}$, the identity (\ref{bar-Psi vs. Psi gen.}) is reduced to 
\begin{equation}
(-1)^{pq}{\bar\Psi}^{i}_{hB{\bar A}j} = \overline{(\Psi^{j}_{A{\bar B}i})}\,. \label{bar-Psi vs. Psi-dagger}
\end{equation}
Unless we explicitly specify the opposite, from now on we will assume that the local frame fields on $E$ are unitary.\\ 

The expression (\ref{bar-Psi vs. Psi-dagger}) motivates the definition of a {\it matricial adjoint} $\Psi^{\dagger}$ of $\Psi$ given by: 
\begin{equation}
\Psi^{\dagger} = \sum\Psi^{\dagger}_{A{\bar B}}\,\theta^{A}\wedge{\bar\theta}^{B} = \sum(\Psi^{\dagger}_{A{\bar B}})^{i}_{j}\,e_{i}\otimes e^{j}\,\theta^{A}\wedge{\bar\theta}^{B}, \label{Psi dagger}
\end{equation}
where $(\Psi^{\dagger}_{A{\bar B}})^{i}_{j} = \overline{(\Psi^{j}_{A{\bar B}i})}$. Hence, using (\ref{bar-Psi vs. Psi-dagger}) and (\ref{Psi dagger}) we have
\begin{equation}
 (-1)^{pq} {\bar\Psi}_{hB{\bar A}} =\Psi^{\dagger}_{A{\bar B}}\,. \label{Key form Psi-bar-h adj.}
\end{equation}
Now, from (\ref{adjoint *}) we have
\begin{eqnarray*}
{\bar *}_{h}\,\Psi &=& *\,\sum{\bar\Psi}_{hB{\bar A}}\,\theta^{B}\wedge\theta^{\bar A}\\
                          &=& \sum{\bar\Psi}_{hB{\bar A}}*(\theta^{B}\wedge\theta^{\bar A})\\
                          &=&  (-1)^{pq}\sum{\bar\Psi}_{hB{\bar A}}\,{\bar *}\,(\theta^{A}\wedge\theta^{\bar B})\,,
\end{eqnarray*}
and consequently from (\ref{Key form Psi-bar-h adj.}) we conclude that
\begin{equation}
{\bar *}_{h}\,\Psi  = \sum\Psi^{\dagger}_{A{\bar B}}\,{\bar *}\,(\theta^{A}\wedge\theta^{\bar B})\,. \label{h.c.Psi}
\end{equation}
At this point we can evaluate the trace of $\Phi\wedge\overline{*}_{h}\Psi$ in a simple way. In fact, a similar procedure to that one developed in Section \ref{Hodge op.} for forms in $A^{p,q}_{X}$, but this time using (\ref{Phi comp.}), (\ref{h.c.Psi}) and (\ref{prop. bar * op. 2}), yields
\begin{eqnarray*}
{\rm tr}(\Phi\wedge{\bar *}_{h}\,\Psi) &=& \sum {\rm tr}\left\{\Phi_{A{\bar B}}\Psi^{\dagger}_{C{\bar D}}\right\}(\theta^{A}\wedge{\bar\theta}^{B})\wedge{\bar *}\,(\theta^{C}\wedge{\bar\theta}^{D})\\
                                                         &=& \sum {\rm tr}\left\{\Phi_{A{\bar B}}\Psi^{\dagger}_{A{\bar B}}\right\}(\theta^{A}\wedge{\bar\theta}^{B})\wedge{\bar *}\,(\theta^{A}\wedge{\bar\theta}^{B})\\
                                                         &=&  \sum {\rm tr}\left\{\Phi_{A{\bar B}}\Psi^{\dagger}_{A{\bar B}}\right\}\frac{\omega^{n}}{n!}\,.                                                       
\end{eqnarray*}
As in the case of forms the above expression implies that 
\begin{equation}
\overline{{\rm tr}(\Phi\wedge{\bar *}_{h}\,\Psi)} = {\rm tr}(\Psi\wedge{\bar *}_{h}\,\Phi)\,, \label{c.c. tr formula}
\end{equation}
and we have a local hermitian inner product $(\cdot\,,\cdot)$ on $A^{p,q}_{X}({\rm End}E)$ given by 
\begin{equation}
(\Phi\,,\Psi) = \sum {\rm tr}\left\{\Phi_{A{\bar B}}\Psi^{\dagger}_{A{\bar B}}\right\}. \label{Local formula Ap,q(EndE)}
\end{equation}
As in (\ref{bar-psi-contrav. comp.}) we define 
\begin{equation}
 {\bar\Psi}_{h}^{\gamma_{1}...\gamma_{p}{\bar\delta}_{1}...{\bar\delta}_{q}} = g^{\beta_{1}{\bar\delta}_{1}}\cdots g^{\beta_{q}{\bar\delta}_{q}} g^{\gamma_{1}{\bar\alpha}_{1}}\cdots g^{\gamma_{p}{\bar\alpha}_{p}}{\bar\Psi}_{h\,\beta_{1}...\beta_{q}{\bar\alpha}_{1}...\bar\alpha_{p}}  \label{bar-Psi-contrav. comp.}
\end{equation}
and if we consider $g^{\alpha\bar\beta}=\delta^{\alpha\bar\beta}$, then (\ref{bar-Psi-contrav. comp.}) becomes
\begin{equation}
{\bar\Psi}_{h}^{A{\bar B}}={\bar\Psi}_{h B{\bar A}}\,. \label{bar-Psi-contrav.} 
\end{equation}
Therefore, from (\ref{Key form Psi-bar-h adj.}) and (\ref{bar-Psi-contrav.}) we can rewrite the local hermitian inner product given by (\ref{Local formula Ap,q(EndE)}) also as:
\begin{equation}
(\Phi\,,\Psi) = (-1)^{pq}\sum{\rm tr}\left\{\Phi_{A{\bar B}}\,{\bar\Psi}_{h}^{A{\bar B}}\right\}. \label{Inner pto. A(EndE) E. Conv.}
\end{equation}
From the above analysis we obtain the following result in complex geometry \cite{Kobayashi}, which is commonly used in any geometric approach to Yang-Mills theory. 
\begin{pro}\label{Proposition A(EndE)}
Let $E$ be a holomorphic vector bundle over a compact K\"ahler manifold $X$ with K\"ahler form $\omega$ and let $h$ be an hermitian metric in $E$, then 
\begin{equation}
\langle\Phi\,,\Psi\rangle = \int_{X}{\rm tr}(\Phi\wedge {\bar *}_{h}\,\Psi) = \int_{X}(\Phi\,,\Psi)\,\frac{\omega^{n}}{n!}\,, \quad\quad \forall \Phi,\Psi\in A^{p,q}_{X}({\rm End}E) \label{Inner pto. A(EndE)}
\end{equation}
gives a global hermitian inner product on $A^{p,q}_{X}({\rm End}E)$, where $(\Phi\,,\Psi)$ is locally given by (\ref{Local formula Ap,q(EndE)}) or equivalently by (\ref{Inner pto. A(EndE) E. Conv.}). 
\end{pro}
As in the case of forms we have local and global norms given by
\begin{equation}
|\Psi|^{2}=(\Psi\,,\Psi)\,,\quad\quad\quad \|\Psi\|^{2}={\langle\Psi\,,\Psi\rangle}\,. \label{2 norms}
\end{equation}
Notice that, in contrast to the theory of Yang-Mills on Riemannian manifolds, in the complex case it is not necessary to include a minus sign in the right hand side of (\ref{Inner pto. A(EndE)}); this is indeed a consequence of the definition of the operator ${\bar *}_{h}$, which carries inside a conjugate transpose operation.\footnote{Notice that if $M, N\in {\rm M}_{n}({\mathbb C})$, the expression ${\rm tr}(MN^{\dagger})$ gives an hermitian inner product in ${\rm M}_{n}({\mathbb C})$. In particular we have ${\rm tr}(MM^{\dagger})=\sum M^{i}_{\,j}(M^{\dagger})^{j}_{\,i} = \sum M^{i}_{\,j}\overline{M^{i}_{j}}\ge 0$.} It is important to note that the hermitian inner products $(\cdot\,,\cdot)$ and $\langle\cdot\,,\cdot\rangle$ depend on the hermitian metric $h$ of $E$ as well as the K\"ahler metric $g$ of $X$. Clearly, as in the case of forms in $A^{p,q}_{X}$, the hermitian inner product given by (\ref{Inner pto. A(EndE)}) can be extended in an obvious way to the space of all forms with coefficients in ${\rm End}E$.\\

The formulas (\ref{Local formula Ap,q(EndE)}) and (\ref{Inner pto. A(EndE) E. Conv.}) are the natural extensions to $A^{p,q}_{X}({\rm End}E)$ of the local formulas (\ref{Local formula Ap,q}) and (\ref{Inner pto. forms E. Conv.}) for $A^{p,q}_{X}$ in complex geometry. Notice that even when (\ref{Inner pto. A(EndE) E. Conv.}) is written following an standard notation in physics, it is eventually more convenient to use (\ref{Local formula Ap,q(EndE)}), e.g., the hermiticity property of the trace (\ref{c.c. tr formula}) is less evident if we use the formula involving the hermitian conjugate instead of the one with the adjoint.


\section{Higgs bundles and the Hermite-Yang-Mills equations}\label{Higgs-HYM section}

In this section we review the basic notions on Higgs bundles and we define the Hermite-Yang-Mills equations. In particular, we give some examples of Higgs bundles. Additionally, we define the Kobayashi and the full Yang-Mills-Higgs functionals for Higgs bundles and we show that there exists a non trivial expression relating these functionals. More details on these topics can be found in \cite{Bruzzo-Granha, Cardona 1, Cardona 2, Simpson} and \cite{Wijnholt}. \\ 

Following Simpson \cite{Simpson}, we define a {\it Higgs bundle} ${\mathfrak E}$ over a compact K\"ahler manifold $X$, as a pair ${\mathfrak E}=(E,\Phi)$, where $E\longrightarrow X$ is a holomorphic vector bundle and $\Phi\in A^{1,0}_{X}({\rm End}E)$ is holomorphic and satisfies the condition: 
\begin{equation}
\Phi\wedge\Phi=0\,. \label{Phi-cond.}
\end{equation}
The ${\rm End}E$-valued form $\Phi$ is commonly called the {\it Higgs field} of the Higgs bundle ${\mathfrak E}$. The first examples of these bundles are considerably technical, e.g., in the article of Hitchin \cite{Hitchin}, such objects appear in the form of bundles associated to the square roots of the canonical bundle of a compact Riemann surface. In the article of Simpson \cite{Simpson}, he defines first the systems of Hodge bundles as objects closely related to the notion of variations of Hodge structures in algebraic geometry. A system of Hodge bundles turns out to be an interesting and non trivial example of a Higgs bundle. To be precise, Simpson defines a system of Hodge bundles as a direct sum of holomorphic bundles $E^{p,q}$ together with maps $\Phi^{p,q}:E^{p,q}\longrightarrow E^{p-1,q+1}$ satisfying 
\begin{equation}
\Phi^{p-1,q+1}\wedge\Phi^{p,q} = 0\,. \nonumber
\end{equation} 
Notice that if $E = \bigoplus E^{p,q}$, the morphisms $\Phi^{p,q}$ define (in the obvious way) a morphism $\Phi$ on $E$ satisfying (\ref{Phi-cond.}). Hence the pair $(E,\Phi)$ becomes a Higgs bundle. \\

Higgs bundles appear also in a natural way from bundles associated to the cotangent bundle of certain compact K\"ahler manifolds $X$ \cite{Seaman}. In fact, suppose that there exists a nowhere vanishing holomorphic form $\lambda\in \Omega^{s,0}_{X}$ with $s$ odd. Notice that such forms can be defined on any Calabi-Yau manifold. Let us consider the holomorphic bundle $E=\bigoplus_{p=0}^{n} \Omega^{p,0}_{X}$ and $\Phi\in A^{1,0}_{X}({\rm End\,}E)$ defined by the condition $\Phi(v)\xi=(\iota_{v}\lambda)\wedge\xi$, where $v$ and $\xi$ are a holomorphic vector field on $X$ and a holomorphic section of $E$. Here 
$\iota_{v}:\Omega^{p,0}_{X}\longrightarrow \Omega^{p-1,0}_{X}$ is the usual contraction operator. Then, $\Phi$ is holomorphic and since $s$ is odd, a direct computation shows that it also satisfies (\ref{Phi-cond.}). Consequently, the pair $(E,\Phi)$ is a Higgs bundle. Moreover, the bundles $E^{a}=\bigoplus_{p\ge a}\Omega^{p,0}_{X}$ with $a\ge 0$ together with morphisms $\Phi$ defined as above furnishes again examples of Higgs bundles. Indeed, such bundles define a filtration of $E$ by Higgs bundles $E^{a}$. For more details on these geometric issues the reader can see \cite{Seaman}.\\  
   
Using the conventions of the preceding sections the general expression (\ref{Phi comp 2.}) gives 
\begin{equation}
\Phi = \sum\Phi_{\alpha}\,\theta^{\alpha} = \sum\Phi^{i}_{\alpha j}\, e_{i}\otimes e^{j}\theta^{\alpha}\,,  \label{Phi components}
\end{equation}
and (\ref{Phi-cond.}) means that 
\begin{equation}
0 = \sum\Phi_{\alpha}\Phi_{\beta}\,\theta^{\alpha}\wedge\theta^{\beta} = \frac{1}{2}\sum_{\alpha<\beta}[\Phi_{\alpha},\Phi_{\beta}]\,\theta^{\alpha}\wedge\theta^{\beta} \nonumber
\end{equation}
and hence (\ref{Phi-cond.}) is equivalent to
\begin{equation}
[\Phi_{\alpha},\Phi_{\beta}] = 0\,. \label{Phi-cond. in components} 
\end{equation}
If $\Phi$ is a Higgs field, its formal adjoint is the element ${\bar\Phi}_{h}$ in $A^{0,1}_{X}({\rm End}E)$ given by (\ref{adjoint Psi}). Therefore 
\begin{equation}
{\bar\Phi}_{h} = \sum {\bar\Phi}_{h{\bar\beta}}\,{\bar\theta}^{\beta} = \sum{\bar\Phi}^{i}_{h{\bar\beta}j}\,e_{i}\otimes e^{j}{\bar\theta}^{\beta}\, \label{barPhi components}
\end{equation}
and the condition (\ref{Phi-cond.}) implies  
\begin{equation}
{\bar\Phi}_{h}\wedge{\bar\Phi}_{h}=0\,, \label{Phi*-cond.}
\end{equation}
which is indeed equivalent to a set of commutation relations for the corresponding endomorphisms ${\bar\Phi}_{h\bar\beta}$. Notice that (\ref{Key form Psi-bar-h adj.}) implies 
\begin{equation}
{\bar\Phi}_{h\bar\beta} = (-1)^{1\cdot0}\Phi^{\dagger}_{\beta}=\Phi^{\dagger}_{\beta} \label{dagger vs. h.c. phi}
\end{equation}
and hence (using unitary frame fields) the matrix representing ${\bar\Phi}_{h\bar\beta}$ is formally the adjoint of the matrix representing $\Phi_{\beta}$.\\

Using the Chern connection $D_{h}$, the Higgs field $\Phi$ and its adjoint ${\bar\Phi}_{h}$, Simpson defines in \cite{Simpson} the connection 
\begin{equation}
{\cal D}_{h} = D_{h} + \Phi + {\bar\Phi}_{h} \label{HS-conn.}
\end{equation} 
which is usually called the {\it Hitchin-Simpson connection}. The curvature of this connection is given by ${\cal F}_{h}={\cal D}_{h}\wedge{\cal D}_{h}$ and is called the 
{\it Hitchin-Simpson curvature}. Using (\ref{Phi-cond.}), (\ref{Phi*-cond.}) and (\ref{HS-conn.}) we have
\begin{eqnarray*}
{\cal F}_{h}  &=&  (D_{h} + \Phi + {\bar\Phi}_{h})\wedge (D_{h} + \Phi + {\bar\Phi}_{h}) \\
                    &=& D_{h}\wedge D_{h} + D_{h}\wedge\Phi + \Phi\wedge D_{h} +  D_{h}\wedge{\bar\Phi}_{h} + {\bar\Phi}_{h}\wedge D_{h}\\
                    && + \; \Phi\wedge{\bar\Phi}_{h} + {\bar\Phi}_{h}\wedge\Phi\\
                    &=& F_{h} + D_{h}\Phi + D_{h}{\bar\Phi}_{h} + [\Phi\,,{\bar\Phi}_{h}]\,.
\end{eqnarray*}
Here $F_{h}=D_{h}\wedge D_{h}$ is the Chern curvature defined in Section \ref{CGeo-HYM}, $D_{h}\Phi$ and $D_{h}{\bar\Phi}_{h}$ are the covariant derivatives of the Higgs field and its hermitian conjugate and we have also used the general commutator formula (\ref{commutator}). Now, since $\Phi$ and ${\bar\Phi}_{h}$ are holomorphic and anti-holomorphic forms respectively, $D_{h}\Phi=D'_{h}\Phi$ and $D_{h}{\bar\Phi}_{h}=d''{\bar\Phi}_{h}$ and the above expression can be further simplified to:
\begin{equation}
{\cal F}_{h}  = F_{h} + D'_{h}\Phi + d''{\bar\Phi}_{h} + [\Phi\,,{\bar\Phi}_{h}]\,. \label{HS-curv.}
\end{equation}  
It is important to note that the curvature given by (\ref{HS-curv.}) has components of different type; to be more precise $D'_{h}\Phi\in A^{2,0}_{X}({\rm End}E)$ and $d''{\bar\Phi}_{h}\in A^{0,2}_{X}({\rm End}E)$ and the remaining part is a form in $A^{1,1}_{X}({\rm End}E)$ given by
\begin{equation}
{\cal F}^{1,1}_{h} = F_{h} + [\Phi\,,{\bar\Phi}_{h}]\,. \label{F11}
\end{equation}

Locally (\ref{F11}) can be written as
\begin{equation}
{\cal F}^{1,1}_{h} = \sum{\cal F}_{h\alpha\bar\beta}\,\theta^{\alpha}\wedge {\bar\theta}^{\beta} =  \sum{\cal F}^{i}_{h\alpha\bar\beta j}\, e_{i}\otimes e^{j} \theta^{\alpha}\wedge{\bar\theta}^{\beta}, \label{cal F components}
\end{equation}
where
\begin{equation}
{\cal F}_{h\alpha\bar\beta} = F_{h\alpha\bar\beta} + [\Phi_{\alpha},{\bar\Phi}_{h\bar\beta}]\,. 
\end{equation}

At this point we can define the {\it Hitchin-Simpson mean curvature} ${\cal K}_{h}$ (see \cite{Cardona 7} for details) as the element in $\Gamma({\rm End}E)=A^{0,0}_{X}({\rm End}E)$ satisfying 
\begin{equation}
{\cal K}_{h}\,\omega^{n} = in\,{\cal F}_{h}\wedge{\omega}^{n-1} =  in\,{\cal F}^{1,1}_{h}\wedge{\omega}^{n-1}. \label{K-R formula}
\end{equation}
Using components we have
\begin{equation}
{\cal K}_{h} = \sum{\cal K}^{i}_{hj}\, e_{i}\otimes e^{j}\,, \quad\quad {\cal K}^{i}_{hj} = \sum g^{\alpha\bar\beta}{\cal F}^{i}_{h\alpha\bar\beta j}\,. 
\end{equation}
Now, in complex geometry it is usual to consider ${\cal K}_{h}$ as an hermitian form $\hat{\cal K}_{h}$ by defining 
\begin{equation}
\hat{\cal K}_{h}(s,s') = h({\cal K}_{h}s,s')\,, \quad\quad \forall s,s'\in\Gamma(E)\,. \label{hat K}
\end{equation} 
An hermitian metric on ${\mathfrak E}$ is said to be {\it Hermite-Yang-Mills (HYM)} or {\it Hermite-Einstein (HE)} \cite{Simpson, Simpson 2} if it satisfies the equation
\begin{equation}
{\cal K}_{h} = c\,I\,, \quad\quad {\rm or \;\, equivalently } \quad\quad \hat{\cal K}_{h} = c\,h\,, \label{HYM-eq.}
\end{equation}
with $c$ the constant given by
\begin{equation}
c=\frac{2\pi\,{\rm deg}E}{r(n-1)!{\rm vol}X}\,.\label{def. c}
\end{equation}
Here ${\rm vol\,}X$ and ${\rm deg\,}E$ are given by (\ref{vol X}) and (\ref{deg E}), respectively. If $\{e_{i}\}_{i=1}^{r}$ is a local frame field on $E$ (not necessarily unitary) we can write (\ref{HYM-eq.}) as 
\begin{equation}
{\cal K}^{i}_{hj} = c\,\delta^{i}_{j}\,, \quad\quad {\rm or \;\, equivalently} \quad\quad \hat{\cal K}_{hi{\bar j}} = c\,h_{i{\bar j}}\,. \label{HYM-eq. local}
\end{equation}
The value of $c$ is indeed a geometric requirement. In fact, from the local expression (\ref{Local commutator}) we have 
\begin{equation}
{\rm tr}[\Phi\,,\bar\Phi_{h}] = \sum{\rm tr}[\Phi_{\alpha},{\bar\Phi}_{h\bar\beta}]\,\theta^{\alpha}\wedge{\bar\theta}^{\beta} = 0\,. \label{tr conm.} 
\end{equation}
Then, taking the trace of (\ref{K-R formula}) and using (\ref{F11}), (\ref{HYM-eq.}) and (\ref{tr conm.}) we get 
\begin{equation}
cr\,\omega^{n} = in\,{\rm tr\,}F_{h}\wedge\omega^{n-1} = 2\pi n\,c_{1}(E)\wedge\omega^{n-1}\,, \label{rel. A}
\end{equation}
where in the last equality we have used the identification $2\pi\,c_{1}(E)=i\,{\rm tr}F_{h}$ (see Section \ref{CGeo-HYM}). At this point, integrating (\ref{rel. A}) and using (\ref{vol X}) and (\ref{deg E}) we get (\ref{def. c}).\\

As far as holomorphic vector bundles is concern, the existence of HYM metrics is equivalent to the notion of Mumford stability,\footnote{It is important to mention that Mumford stability is just a particular kind of stability defined for holomorphic vector bundles. In fact, in complex geometry there exist other notions of stability, e.g., T-stability and Gieseker stability \cite{Kobayashi} and these stabilities can be also extended to the Higgs bundle case \cite{Cardona 4, Cardona-Mata}.} such equivalence is a remarkable fact in complex geometry which is commonly known as the Hitchin-Kobayashi correspondence. In particular, Simpson proves in \cite{Simpson} an extension of the Hitchin-Kobayashi correspondence for Higgs bundles. In physics literature this result is also known as the Donaldson-Uhlenbeck-Yau theorem \cite{Wijnholt}. To be precise, the result of Simpson establishes that a Higgs bundle has an HYM metric if and only if it is Mumford poly-stable \cite{Simpson}. We are not going to address these algebraic aspects here, for more details on this part the reader can see the pioneering articles of Simpson \cite{Simpson, Simpson 2}.\\

As in the classical case of holomorphic vector bundles, there exist some functionals of interest that can be defined in the space of hermitian metrics of Higgs bundles. In particular, Simpson defines a Donaldson functional for Higgs bundles in \cite{Simpson}. The Donaldson functional can be also introduced following the geometric approach of Kobayashi \cite{Kobayashi}. It is important to mention that such a functional plays an crucial role in the proof of the Hitchin-Kobayashi correspondence for Higgs bundles (see \cite{Cardona 1, Cardona 2} for details). In fact, an hermitian metric in a Higgs bundle is a critical point of the Donaldson functional if and only if it is HYM. More details on the Donaldson functional and the Hitchin-Kobayashi correspondence for Higgs bundles can be found in \cite{Cardona 1, Cardona 2, Simpson}.\\

On the other hand, a natural functional that can be defined for Higgs bundles is the {\it full Yang-Mills-Higgs functional}. This is defined as the square of the norm of the Hitchin-Simpson curvature ${\cal F}_{h}$, where the norm is computed from the global hermitian inner product of Proposition \ref{Proposition A(EndE)}. Hence, using (\ref{HS-curv.}) this functional is:
\begin{equation}
\|{\cal F}_{h}\|^{2} = \| F_{h} + [\Phi\,,{\bar\Phi}_{h}] \|^{2} + \|D'_{h}\Phi\|^{2} + \|d''{\bar\Phi}_{h}\|^{2}\,,\label{YMH fnal.}
\end{equation}
or equivalently 
\begin{equation}
\|{\cal F}_{h}\|^{2} = \int_{X}\left\{|F_{h}+[\Phi\,,{\bar\Phi}_{h}]|^{2} + |D'_{h}\Phi |^{2} + |d''{\bar\Phi}_{h}|^{2} \right\}\frac{\omega^{n}}{n!}\,. \label{YMH fnal. int.}
\end{equation}
 
Notice that if $\Phi\equiv 0$, then $\|{\cal F}_{h}\|^{2}=\|F_{h}\|^{2}$ and the full Yang-Mills-Higgs functional becomes the usual Yang-Mills functional for the corresponding holomorphic vector bundle $E\longrightarrow X$ (see \cite{Kobayashi}, Ch. 4 for details).\\ 

Another functional of interest in the case of Higgs bundles is the {\it functional of Kobayashi}, which is indeed proportional to the square of the norm of the Hitchin-Simpson mean curvature ${\cal K}_{h}$. More precisely
\begin{equation}
{\cal J}(h) = \frac{n!}{2}\|{\cal K}_{h}\|^{2} = \frac{1}{2}\int_{X}|{\cal K}_{h}|^{2}\omega^{n}\,. \label{J fnal.} 
\end{equation}  
The functional of Kobayashi satisfies some interesting properties. Clearly, from the definition it is non-negative and it can be shown that for any hermitian metric $h$ we have 
\begin{equation}
{\cal J}(h) \ge \frac{2n(\pi{\rm deg\,}E)^{2}}{r(n-1)!{\rm vol\,}X} \label{J prop.}
\end{equation}   
and that ${\cal J}$ attains this lower bound at $h=h_{0}$ if and only if $h_{0}$ is HYM. It can be shown also that
\begin{eqnarray}
 \| F_{h} + [\Phi\,,{\bar\Phi}_{h}] \|^{2} - \|{\cal K}_{h}\|^{2} &=& 4\pi^{2} \int_{X}\left\{2c_{2}(E) - c_{1}(E)^{2}\right\}\wedge\frac{\omega^{n-2}}{(n-2)!}\nonumber\\
                                                                                              && + \; 2\int_{X}{\rm tr}(F_{h}\wedge[\Phi\,,{\bar\Phi}_{h}])\wedge\frac{\omega^{n-2}}{(n-2)!}. \label{Koba vs I}
\end{eqnarray}
The inequality \eqref{J prop.} and the difference \eqref{Koba vs I} appear in \cite{Cardona 7} as Theorem 1 and Proposition 1, respectively.\footnote{In fact, Proposition 1 in \cite{Cardona 7} involves the trace ${\rm tr}({\cal F}^{1,1}_{h}\wedge [\Phi\,,{\bar\Phi}_{h}])$, but since ${\rm tr}([\Phi\,,{\bar\Phi}_{h}]^{2})=0$, such a trace is equal to ${\rm tr}(F_{h}\wedge[\Phi\,,{\bar\Phi}_{h}])$.} Notice that in the right hand side of \eqref{Koba vs I}, the first term does not depend on $h$ (it is a topological constant), however the second term depends on $h$.\\

In the particular case when $\Phi\equiv 0$, the Higgs bundle is identified with the holomorphic vector bundle $E\longrightarrow X$ and second term in the right hand side of \eqref{Koba vs I} vanishes. Therefore, as far as holomorphic vector bundles is concern (and up to constants) the Yang-Mills and the Kobayashi functionals are the same functional. For Higgs bundles, the situation is more subtle and the corresponding functionals in the above difference are two different functionals. To be precise, from the decomposition (\ref{YMH fnal.}) and the difference \eqref{Koba vs I} we get the following result.

\begin{pro}\label{Prop. FullYM vs Koba}
Let ${\mathfrak E}=(E,\Phi)$ be a Higgs bundle over a compact K\"ahler manifold $X$ with K\"ahler form $\omega$. Let $h$ be an hermitian metric in ${\mathfrak E}$ with ${\cal F}_{h}$ and ${\cal K}_{h}$ the corresponding Hitchin-Simpson curvature and mean curvature of ${\mathfrak E}$, respectively. Then 
\begin{eqnarray}
 \| {\cal F}_{h}\|^{2} - \|{\cal K}_{h}\|^{2} &=& \|D'_{h}\Phi\|^{2} + \|d''{\bar\Phi}_{h}\|^{2} + 2\int_{X}{\rm tr}(F_{h}\wedge[\Phi\,,{\bar\Phi}_{h}])\wedge\frac{\omega^{n-2}}{(n-2)!}\nonumber\\
                                                               &&+\; 4\pi^{2}\int_{X}\left\{2c_{2}(E) - c_{1}(E)^{2}\right\}\wedge\frac{\omega^{n-2}}{(n-2)!}\,. \label{Koba vs fullYM}
\end{eqnarray}
\end{pro}
The formula \eqref{Koba vs fullYM} in Proposition \ref{Prop. FullYM vs Koba} is an expression relating the full Yang-Mills-Higgs and Kobayashi functionals for Higgs bundles in a non trivial way. In fact, the right hand side of \eqref{Koba vs fullYM} depends on three terms involving $h$. In other words, if we consider the hermitian metric as a variable, the difference between the full Yang-Mills-Higgs and the Kobayashi functionals is not constant, but depends on the hermitian metric $h$. 


\section{Higgs bundles and $2k$-Hitchin's equations}\label{Higgs-Ward section}

The main purpose of this section is to study the $2k$-Hitchin equations introduced by Ward \cite{Ward} from the geometric viewpoint of Higgs bundles. In particular, we show that as far as Higgs bundles is concern, the $2k$-Hitchin equations are reduced to a set of two equations defined in the space of hermitian metrics of the Higgs bundle. Finally, we show that if there exists an hermitian metric satisfying the $2k$-Hitchin equations, then it is a minimum of certain funcional. \\ 

From a geometric viewpoint \cite{Lubke}, we can fix the holomorphic structure on a complex vector bundle $E$ and consider all possible hermitian structures on it, i.e., we can consider all hermitian metrics in $E$; or we can do the opposite, namely, we can fix the hermitian metric $h$ in $E$ and consider all possible holomorphic structures on $E$. In this article we will consider $2k$-Hitchin's equations following the geometric approach of Kobayashi \cite{Kobayashi}. Therefore, we begin by fixing a holomorphic structure on a complex vector bundle and hence we can think such a bundle as a holomorphic vector bundle $E\longrightarrow X$, where $X$ is a compact K\"ahler manifold. Then, we consider all possible hermitian metrics $h$ in $E$ and hence some  objects will depend on $h$. From this perspective $2k$-Hitchin's equations (\ref{2k-Hitchin eqs.}) will be rewritten as:
\begin{equation}
D_{h}\Phi = 0\,, \quad\quad F^{1,1}_{h} + \frac{1}{4}[\Phi\,,{\bar\Phi}_{h}] = 0\,, \quad\quad [\Phi\,,\Phi] = 0\,, \quad\quad F^{2,0}_{h} = 0\,. \label{2k-Hitchin eqs. rev.}
\end{equation} 
From this point of view, the equations are defined for pairs $(h,\Phi)$ with $h$ an hermitian metric in $E$ and $\Phi\in A^{1,0}_{X}({\rm End}E)$.\\

Notice that if $D_{h}$ is the Chern connection of $E$ defined by $h$, then its Chern curvature $F_{h}$ is an element in $A^{1,1}_{X}({\rm End\,}E)$. Then $F_{h}=F^{1,1}_{h}$ and $F^{2,0}_{h}=0$ and the last equation in (\ref{2k-Hitchin eqs. rev.}) is satisfied. Now, the third equation is equivalent to $\Phi\wedge\Phi=0$, and hence it will be also satisfied if $\Phi$ is considered as a Higgs field of $E$. Indeed, such a condition is imposed on any Higgs field in higher dimensions \cite{Simpson}. In summary, if we consider the Chern connection $D_{h}$ and a Higgs field $\Phi$ of $E$, the last two equations in (\ref{2k-Hitchin eqs. rev.}) will be ``formally satisfied''. Since $\Phi$ is now a Higgs field, it is necessarily holomorphic and $d''\Phi=0$. Then the first equation can be further reduced to $D'_{h}\Phi=0$. Hence (after a rescaling of the Higgs field) the $2k$-Hitchin equations (\ref{2k-Hitchin eqs. rev.}) are reduced to:
\begin{equation}
D'_{h}\Phi = 0\,, \quad\quad F_{h} + [\Phi\,,{\bar\Phi}_{h}] = 0\,, \label{2k-Hitchin eqs. revisited}
\end{equation}
where ${\bar\Phi}_{h}\in A^{0,1}_{X}({\rm End}E)$ is the hermitian conjugate of $\Phi$ with respect to $h$ defined in (\ref{adjoint Psi}). From a geometric viewpoint, the first equation in (\ref{2k-Hitchin eqs. revisited}) gives a parallelism condition on $\Phi\in A^{1,0}_{X}({\rm End}E)$ with respect to the Chern connection $D_{h}$, and the second equation is a constraint of the curvature and the Higgs field. From the above analysis, it is clear that Higgs bundles ${\mathfrak E}=(E,\Phi)$ seem to be a natural setting for studying the $2k$-Hitchin equations. In fact, using such bundles the question is reduced to determine if a Higgs bundle admits or not an hermitian metric $h$ satisfying (\ref{2k-Hitchin eqs. revisited}). Now, as we already mentioned, the existence of HYM metrics on Higgs bundles is closely related to notions of stability on these bundles. Hence, it is natural to wonder if there exists a notion of stability related to the existence of hermitian metrics satisfying the $2k$-Hitchin equations (\ref{2k-Hitchin eqs. revisited}). We do not know at the moment if the existence of solutions to such equations is related to one of the notions of stability studied in algebraic geometry. Such issues are beyond the scope of this survey and we are not going to address these questions here. For more details on Higgs bundles and the role of stability in theoretical physics the reader can see \cite{Wijnholt}.\\
 
We can write (\ref{2k-Hitchin eqs. revisited}) in components. In fact, using a unitary local frame field $\{\theta^{\alpha}\}_{\alpha=1}^{n}$ of $\Omega^{1,0}_{X}$ we can write $\Phi$ and ${\bar\Phi}_{h}$ as in (\ref{Phi components}) and (\ref{barPhi components}) and 
\begin{equation}
F_{h} = \sum F_{h\alpha{\bar\beta}}\,\theta^{\alpha}\wedge{\bar\theta}^{\beta} = \sum F^{i}_{h\alpha{\bar\beta}j}\,e_{i}\otimes e^{j}\,\theta^{\alpha}\wedge{\bar\theta}^{\beta}\,. \label{F comp.}
\end{equation}
Notice that $F_{h}$ is the the Chern curvature, then from (\ref{F11}) it is clear that the second equation in (\ref{2k-Hitchin eqs. revisited}) can be written in terms of the Hitchin-Simpson curvature as ${\cal F}^{1,1}_{h}=0$. If we denote by $A_{h}$ the connection form of $D'_{h}$ the equations (\ref{2k-Hitchin eqs. revisited}) become
\begin{equation}
\partial_{\alpha}\Phi_{\beta}  + [A_{h\alpha},\Phi_{\beta}] = 0\,, \quad\quad\quad  F_{h\alpha{\bar\beta}} + [\Phi_{\alpha}\,,{\bar\Phi}_{h\bar\beta}] = 0\,,  \label{2k-Hitchin eqs. rev. comp.}
\end{equation}
where the indices $\alpha,\beta = 1,...,n$. These can be written even in a more explicit way if we consider a unitary local frame field $\{e_{i}\}_{i=1}^{r}$ of $E$. In fact, using such a frame (\ref{2k-Hitchin eqs. rev. comp.}) are  
\begin{equation}
\partial_{\alpha}\Phi^{i}_{\beta j} + \sum (A^{i}_{h\alpha k}\Phi^{k}_{\beta j} - \Phi^{i}_{\beta k}A^{k}_{h\alpha j})= 0\,, \label{2k-Hitchin eqs. rev. comp. 1}
\end{equation}
\begin{equation}
F^{i}_{h\alpha{\bar\beta}j} + \sum (\Phi^{i}_{\alpha k}{\bar\Phi}^{k}_{h\bar\beta j} -  {\bar\Phi}^{i}_{h\bar\beta k}\Phi^{k}_{\alpha j}) = 0\,, \label{2k-Hitchin eqs. rev. comp. 2}
\end{equation}
where the indices $\alpha,\beta = 1,...,n$ and $i,j = 1,...,r$.\\ 

The equations (\ref{2k-Hitchin eqs. revisited}) are similar to the Seiberg-Witten equations \cite{Friedrich}.\footnote{A brief introduction to Seiberg-Witten theory can be found in the Appendix A of such a reference.} In fact, in both cases we have a connection form and a certain field, and the 
equations are given by a parallelism condition of the field with respect to the covariant derivative of the connection form and a constraint equation involving the curvature of the covariant derivative and the field. The main difference between the Seiberg-Witten equations and $2k$-Hitchin's equations comes from the ``nature'' of the field; meanwhile in the former equations it is an spinor field, in the later equations it is a Higgs field. Following the analogy with the Seiberg-Witten theory, we associate to (\ref{2k-Hitchin eqs. revisited}) the functional
\begin{equation}
{\cal H}(h) = \| D'_{h}\Phi \|^{2} + \| F_{h} + [\Phi\,,{\bar\Phi}_{h}]\|^{2}\,. \label{cal H}
\end{equation}
We call the functional ${\cal H}$ the {\it non-abelian Seiberg-Witten fucntional}, such a functional is defined in the space ${\rm Herm}^{+}{\mathfrak E}$ of hermitian metrics of ${\mathfrak E}$. This space is formally the same than the space of the hermitian metrics of the holomorphic vector bundle $E$ \cite{Cardona 1, Cardona 2} (see \cite{Kobayashi}, Ch. VI for basic definitions and properties of this space). At this point, using the local norm $|\cdot|$ introduced in Section \ref{CGeo-HYM}, the functional (\ref{cal H}) can be written as
\begin{equation}
{\cal H}(h) = \int_{X}\left\{\lvert D'_{h}\Phi\lvert^{2} + \lvert F_{h} + [\Phi\,,{\bar\Phi}_{h}]\lvert^{2}\right\}\frac{\omega^{n}}{n!}\,.  \label{cal H omega}
\end{equation}
Now, notice that from (\ref{YMH fnal.}) and (\ref{cal H}) we get
\begin{equation}
\|{\cal F}_{h}\|^{2} = {\cal H}(h) +  \|d''{\bar\Phi}_{h}\|^{2},
\end{equation} 
which shows that the full Yang-Mills-Higgs functional and the non-abelian Seiberg-Witten functional ${\cal H}$ are closely related. However, this relation is far away to be trivial, since the difference between these functionals is the global norm of the term $d''{\bar\Phi}_{h}$, which depends on $h$.\\
 
On the other hand, using (\ref{cal H}) or (\ref{cal H omega}) it is clear that ${\cal H}$ is non-negative and just by construction of the functional we have the following result.
\begin{pro}\label{Prop. 2k-sol. minimum of H}
Let ${\mathfrak E}=(E,\Phi)$ be a Higgs bundle over a compact K\"ahler manifold $X$. If there exists an hermitian metric $h_{0}$ in $E$ satisfying the $2k$-Hitchin equations (\ref{2k-Hitchin eqs. revisited}), then $h_{0}$ is a minimum of the functional ${\cal H}$ given by (\ref{cal H}).
\end{pro}
From a geometric viewpoint, for a Higgs bundle it is important to know not only the minima of a functional but also the critical points of such functional \cite{Cardona 7}. Therefore, it will be useful to apply some variational technics to ${\cal H}$ in order to find the Euler-Lagrange equations associated with the functional ${\cal H}$. This study will be beyond the scope of this survey and we hope to come back to such questions in a forthcoming article.\\

Finally, the terms in the integrand of (\ref{cal H omega}) can be easily computed using the theory developed in Section \ref{CGeo-HYM}. In fact, from 
(\ref{Key form Psi-bar-h adj.}) it follows that $(\partial_{\alpha}{\Phi}_{\beta})^{\dagger}=\partial_{\bar\alpha}{\bar\Phi}_{h\bar\beta}$ and hence using (\ref{Local formula Ap,q(EndE)}), (\ref{bar-Psi-contrav.}) and (\ref{dagger vs. h.c. phi}) we get
\begin{eqnarray}
\lvert D'_{h}\Phi \lvert ^{2} &=& \sum{\rm tr} \left\{\left(\partial_{\alpha}\Phi_{\beta} + [A_{h\alpha},\Phi_{\beta}]\right)(\partial_{\alpha}{\Phi}_{\beta} +  [A_{h\alpha},\Phi_{\beta}])^{\dagger}\right\} \nonumber \\
                                          &=& \sum{\rm tr} \left\{\left(\partial_{\alpha}\Phi_{\beta} + [A_{h\alpha},\Phi_{\beta}]\right)(\partial^{\alpha}{\bar\Phi}_{h}^{\beta} - [{\bar A}_{h}^{\alpha},{\bar\Phi}_{h}^{\beta}])\right\}, \label{Norm a}
\end{eqnarray} 
where ${\bar A}_{h}$ denotes the hermitian conjugate of $A_{h}$. Now, since $(\cdot\,,\cdot)$ is a local hermitian inner product we have
\begin{eqnarray*}
   \lvert F_{h} + [\Phi\,,{\bar\Phi}_{h}]\lvert ^{2} &=& ( F_{h} + [\Phi\,,{\bar\Phi}_{h}]\,,  F_{h} + [\Phi\,,{\bar\Phi}_{h}] )\\
                                                                       &=& \lvert F_{h}\lvert ^{2} + 2\,{\mathfrak{Re}}\,([\Phi\,,{\bar\Phi}_{h}]\,,F_{h}) + \lvert [\Phi\,,{\bar\Phi}_{h}] \lvert ^{2}\,, 
\end{eqnarray*}
where ${\mathfrak{Re}}$ is the real part of the corresponding inner product. At this point, we can write the three terms in the above expression in a similar way as we have done in (\ref{Norm a}). After doing that the functional (\ref{cal H omega}) becomes
\begin{equation}
{\cal H}(h) = \int_{X}{\cal L}(A_{h},\Phi,{\bar\Phi}_{h})\,\frac{\omega^{n}}{n!}\,, \label{H loc. 1}
\end{equation}
with 
 \begin{eqnarray*} 
 {\cal L}(A_{h},\Phi,{\bar\Phi}_{h}) &=& \sum{\rm tr} \left\{\left(\partial_{\alpha}\Phi_{\beta} + [A_{h\alpha},\Phi_{\beta}]\right)(\partial^{\alpha}{\bar\Phi}_{h}^{\beta} - [{\bar A}_{h}^{\alpha},{\bar\Phi}_{h}^{\beta}]) - F_{h\alpha\bar\beta}{\bar F}_{h}^{\alpha\bar\beta}\right\}\\
                                                     &-&  \sum{\rm tr}\left\{[\Phi_{\alpha},{\bar\Phi}_{h\bar\beta}]\,[{\bar\Phi}_{h}^{\alpha},\Phi^{\bar\beta}] + 2\,{\mathfrak{Re\,}}([\Phi_{\alpha},{\bar\Phi}_{h\bar\beta}]\,{\bar F}_{h}^{\alpha\bar\beta})\right\}.
\end{eqnarray*} 
In the above expression ${\bar F}_{h}$ denotes the hermitian conjugate of $F_{h}$. From the physical point of view, the lagrangian ${\cal L}(A_{h},\Phi,{\bar\Phi}_{h})$ represents a theory involving interactions between a gauge field $A_{h}$ and a Higgs field $\Phi$ in the same spirit of the lagrangians appearing in the celebrated works of Kapustin and Witten \cite{Witten-Kapustin, Witten 2}.


\section{Appendix: The origin of Hitchin's equations}\label{Appendix}

For the benefit of the reader, in this section we review how the Hitchin equations \cite{Hitchin} arise as a consequence of a dimensional reduction procedure applied to the self-dual Yang Mills equations on ${\mathbb R}^{4}$. Additionally, we present some equivalent forms in which these equations can be written, indeed, some of them are the ways in which Hitchin's equations appear in physics literature \cite{Witten-Kapustin, Ward, Ward 2}.\\

Let us consider the Riemannian manifold ${\mathbb R}^{4}$ with the usual metric and coordinates given by $g_{ij}=\delta_{ij}$ and $x^{i}$ (with $i,j=1,...,4$). Let $A = A_{i}\,dx^{i}$ be a $SU(2)$ (smooth) gauge potential on ${\mathbb R}^{4}$ with corresponding covariant derivatives $D_{i} = \partial_{i} + A_{i}$ and gauge field $F=\frac{1}{2}\sum F_{ij}\,dx^{i}\wedge dx^{j}$. As it is well known, the components of the gauge field are given by 
\begin{equation}
F_{ij} = [D_{i},D_{j}] = \partial_{i}A_{j} -   \partial_{j}A_{i} + [A_{i},A_{j}]\,. \label{comp. F}
\end{equation}
Considering the usual Hodge $*$ operator on ${\mathbb R}^{4}$ associated to $g_{ij}$ we have that $*F$ is again a 2-form.\footnote{If $\epsilon_{i_{1}....i_{4}}$ is the usual Levi-Civita symbol on ${\mathbb R}^{4}$, the $*$ operator in Differential Geometry is given by 
\begin{equation}
*(dx^{i_{1}}\wedge...\wedge dx^{i_{p}}) = \frac{1}{(4-p)!}\sum g^{i_{1}j_{1}}...g^{i_{p}j_{p}}\epsilon_{j_{1}...j_{p}j_{p+1}...j_{4}}\,dx^{j_{p+1}}\wedge....\wedge dx^{j_{4}}.
\end{equation}
In particular we have $*(dx^{1}\wedge dx^{2})=dx^{3}\wedge dx^{4}$ and $*(dx^{1}\wedge dx^{3})=dx^{4}\wedge dx^{2}$ and $*(dx^{1}\wedge dx^{4})=dx^{2}\wedge dx^{3}$.}
A gauge potential is called {\it self-dual Yang-Mills} (SDYM) if $F$ is invariant under the Hodge operator, i.e., $*F = F$ or more explicitly, if in terms of the  components $F_{ij}$ we have
\begin{equation}
F_{12} = F_{34}\,,\quad\quad F_{13} = F_{42}\,,\quad\quad F_{14} = F_{23}\,. \label{SDYM eqs.}
\end{equation} 
If we assume that each $A_{i}$ does not depend on two coordinates, say $A_{i}=A_{i}(x^{1},x^{2})$, we have a dimensional reduction from 
${\mathbb R}^{4}$ to ${\mathbb R}^{2}$, where now $A=A_{1}dx^{1} + A_{2}dx^{2}$ and $F=F_{12}\,dx^{1}\wedge dx^{2}$ are  interpreted as a new gauge potential and field on ${\mathbb R}^{2}$ and $A_{3}=\phi_{1}$ and $A_{4}=\phi_{2}$ are ``auxiliary'' fields usually called {\it Higgs fields} \cite{Hitchin, Ward, Ward 2}. This terminology comes from physics and it is indeed standard in the dimensional reduction procedure; however, it is important to note that these Higgs fields are not directly related to the Higgs boson of the standard model. In fact, strictly speaking and from a physical point of view, the Higgs fields here take values in the adjoint representation of such a gauge group, in contrast with the Higgs boson, which takes values in the fundamental representation of $SU(2)$. By using (\ref{comp. F}) the equations (\ref{SDYM eqs.}) can be rewritten as:
\begin{equation}
[D_{1},D_{2}] = [\phi_{1},\phi_{2}]\,, \quad\quad [D_{1},\phi_{1}] = [\phi_{2},D_{2}]\,, \quad\quad [D_{1},\phi_{2}] = [D_{2},\phi_{1}]\,. 
\label{red. SDYM}
\end{equation}
The equations (\ref{red. SDYM}) are the so called {\it Hitchin's equations} and have played an important role in complex geometry and mathematical physics since the 80's. These can be written in different ways; in fact, since the Higgs fields $\phi_{i}$ with $i=1,2$ take values in ${\mathfrak su}(2)$, these fields are given by traceless anti-hermitian matrices and hence $\phi_{i}^{*}=-\phi_{i}$, where the superscript $*$ represents here the transpose conjugate. At this point, if we define a complex Higgs field 
$\phi = \phi_{1} - i\phi_{2}$ we get $\phi^{*}=-\phi_{1} -i\phi_{2}$ and (\ref{red. SDYM}) becomes
\begin{equation}
F_{12} = \frac{i}{2}[\phi,\phi^{*}]\,,\quad\quad (D_{1} + iD_{2})\phi =0\,, \label{Hitchin real}
\end{equation}
where in the second equation we have $D_{j}\phi= [D_{j},\phi]$. Now, if we introduce the complex variable $z=x^{1} +ix^{2}$, the complex (anti-holomorphic) derivative $\partial_{\bar z} = \partial_{1} + i\partial_{2}$ and $A_{\bar z} = A_{1} + i A_{2}$, we get $D_{1} + iD_{2} = \partial_{\bar z} + A_{\bar z} = D_{\bar z}$ which can be seen as a anti-holomorphic covariant derivative. In terms of this complex variable, the gauge field is $F=F_{z{\bar z}}\,dz\wedge d\bar z$ with $F_{z\bar z}=\frac{i}{2}F_{12}$ and (\ref{Hitchin real}) are finally written as:
\begin{equation}
F_{z\bar z} = -\frac{1}{4}[\phi,\phi^{*}]\,,\quad\quad D_{\bar z}\phi =0\,. \label{Hitchin complex}
\end{equation}
The equations (\ref{Hitchin real}) and (\ref{Hitchin complex}) are two ways of writting the Hitchin equations that have been considered in mathematical physics literature; in particular, these equations appear in recently works by Ward \cite{Ward, Ward 2}. It is important to mention that there is also another way in which the Hitchin equations appear in complex geometry and which is a reformulation of (\ref{Hitchin real}) and (\ref{Hitchin complex}) using complex differential forms. To be more precise, by defining 
$\Phi_{c} = \frac{1}{2}\phi\,dz$ and $\Phi^{*}_{c} = \frac{1}{2}\phi^{*}d\bar z$, we have $\Phi_{c}$ and $\Phi^{*}_{c}$ as $(1,0)$ and $(0,1)$ complex differential forms with coefficients in 
${\mathfrak su}(2)$ and the Hitchin equations (\ref{Hitchin complex}) can be written as:
\begin{equation}
F = -[\Phi_{c},\Phi^{*}_{c}]\,, \quad\quad D_{\bar z}\Phi_{c}=0\,, \label{Hitchin complex forms}
\end{equation}  
where the bracket here is the usual extension of the commutator to matrix valued differential forms. In fact, it is a particular case of the general commutator defined by (\ref{commutator}). Hence, one has 
\begin{equation}
[\Phi_{c},\Phi^{*}_{c}] = \Phi_{c}\wedge\Phi^{*}_{c} + \Phi^{*}_{c}\wedge\Phi_{c}\,.  \nonumber
\end{equation}
The equations (\ref{Hitchin complex forms}) are the way in which the Hitchin equations are written in \cite{Hitchin}. In fact, under the identification $D_{\bar z}=d''_{A}$ and ignoring the subscript $c$ on the Higgs field (\ref{Hitchin complex forms}) is exactly  (\ref{Hitchin eqs. 0}). Is in this form that the equations were generalized latter on by Simpson \cite{Simpson} to a set of equations usually called the Hermite-Yang-Mills equations. These generalized equations play an important role in complex geometry as well as in mathematical physics \cite{Simpson, Simpson 2, Wijnholt}.\\

On the other hand, as it is shown in \cite{Cardona-Compean-Merino}, if instead of considering a complex Higgs field $\Phi_{c}$ we define the forms $\Phi =\phi_{1}dx^{1} + \phi_{2}dx^{2}$, then (\ref{Hitchin complex forms}) can be written as:
\begin{equation}
F - \Phi\wedge\Phi = 0\,, \quad\quad D\Phi =0\,, \quad\quad D^{*}\Phi =0\,,   \label{Hitchin real forms}
\end{equation}
where $D = D_{1}dx^{1} + D_{2}dx^{2}$ is the covariant derivative of $A = A_{1}dx^{1} + A_{2}dx^{2}$ and $D^{*}=*D*$. The Hitchin equations written in the form (\ref{Hitchin real forms}) appear in a celebrated work of Kapustin and Witten \cite{Witten-Kapustin}, in which a physical approach to the geometric Langlands program is proposed.\\

Finally, it is important to note that since $D_{\bar z}\phi = \partial_{\bar z}\phi + [A_{\bar z},\phi]$, the Hitchin equations impose a holomorphic condition on the determinant of 
$\phi$. In fact, from the second equation in (\ref{Hitchin complex}) we get 
\begin{equation}
{\rm det}D_{\bar z}\phi = {\rm det\,}\partial_{\bar z}\phi + {\rm det}[A_{\bar z},\phi]  = \partial_{\bar z}\,{\rm det\,}\phi = 0 
\end{equation}
and hence ${\rm det}\,\phi$ is a holomorphic function in $z$. This is indeed a crucial constraint of the Higgs field used by Ward in \cite{Ward, Ward 2}.


\end{document}